\documentclass{elsart3-1}

\usepackage{amssymb,amsmath}

\usepackage[english,francais]{babel}

\newtheorem{theorem}{Theorem}[section]
\newtheorem{lemma}[theorem]{Lemma}
\newtheorem{e-proposition}[theorem]{Proposition}
\newtheorem{corollary}[theorem]{Corollary}
\newtheorem{e-definition}[theorem]{Definition\rm}

\renewcommand{\theequation}{\arabic{equation}}
\def\R{\mathbb{R}}
\def\RR{\mathbb{R}}
\def\di{\displaystyle}

\def\og{\leavevmode\raise.3ex\hbox{$\scriptscriptstyle\langle\!\langle$~}}
\def\fg{\leavevmode\raise.3ex\hbox{~$\!\scriptscriptstyle\,\rangle\!\rangle$}}
\def\p{\partial}
\def\e{\varepsilon}

\journal{the Acad\'emie des sciences}
\begin{document}
\centerline{}
\begin{frontmatter}


\selectlanguage{english}
\title{Uniqueness in a class of Hamilton-Jacobi equations with constraints}


\selectlanguage{english}
\author[authorlabel1]{Sepideh Mirrahimi},
\ead{sepideh.mirrahimi@math.univ-toulouse.fr}
\author[authorlabel1]{Jean-Michel Roquejoffre}
\ead{jean-michel.roquejoffre@mth.univ-toulouse.fr}

\address[authorlabel1]{Institut de Math\'ematiques (UMR CNRS 5219), Universit\'e Paul Sabatier, 118 Route de Narbonne, 31062 Toulouse Cedex, France}


\medskip
\begin{center}
{\small Received *****; accepted after revision +++++\\
Presented by £££££}
\end{center}

\begin{abstract}
\selectlanguage{english}
In this note, we discuss a class of time-dependent Hamilton-Jacobi equations depending on a function of time, this function being chosen in order to keep the maximum of the solution to the constant value 0. The main result of the note is that the full problem has a unique classical solution. The motivation is a selection-mutation model which, in the limit of small diffusion, exhibits concentration on the zero level set of the solution of the Hamilton-Jacobi equation. The uniqueness result that we prove implies strong convergence and error estimates for the selection-mutation model. {\it To cite this article: S. Mirrahimi, J.-M. Roquejoffre, C. R. Acad. Sci. Paris, Ser. I 340 (2015).}

\vskip 0.5\baselineskip

\selectlanguage{francais}
\noindent{\bf R\'esum\'e} \vskip 0.5\baselineskip \noindent
{\bf Unicit\'e pour une classe d'\'equations de Hamilton-Jacobi avec contraintes. }
Dans cette note, on discute une classe d'\'equations de Hamilton-Jacobi d\'ependant du temps, et d'une fonction inconnue du temps choisie pour que le maximum de la 
solution de l'\'equation de Hamilton-Jacobi prenne tout le temps la valeur 0. Le r\'esultat principal de cette note est que le probl\`eme complet admet une unique solution 
classique. La motivation est un mod\`ele de s\'election-mutation qui, dans la limite d'une diffusivit\'e nulle, présente une concentration sur la ligne de niveau 0 de la solution de l'\'equation de Hamilton-Jacobi. Le résultat d'unicit\'e que nous d\'emontrons implique une convergence forte avec estimations d'erreur pour le mod\`ele de s\'election-mutation.
{\it Pour citer cet article~: S. Mirrahimi, J.-M. Roquejoffre, C. R. Acad. Sci.
Paris, Ser. I 340 (2015).}

\end{abstract}
\end{frontmatter}

\selectlanguage{francais}
\section*{Version fran\c{c}aise abr\'eg\'ee}
On pr\'esente dans cette note un r\'esultat d'unicit\'e pour le probl\`eme de Hamilton-Jacobi  suivant,
d'inconnues $(u(t,x),I(t))$:
\begin{equation}
\label{e0.1}
\left\{
\begin{array}{rll}
\partial_t u &= | \nabla  u |^2 +R(x,I), \  \ (t>0,x\in\R^d), \  \ \max_{x} u(t,x) =0,\\
I(0)&=I_0>0,\  \  u(0,x)=u_0(x).
\end{array}
\right.
\end{equation}
La donn\'ee $R(x,I)$ v\'erifie des hypothèses de stricte concavit\'e par rapport \`a $x$ et de monotonie par rapport \`a $I$ explicitées plus bas. La donn\'ee initiale $u_0$ v\'erifie aussi des hypoth\`eses sp\'eciales de concavit\'e. Ainsi, la fonction $I(t)$ doit \^etre choisie pour que la solution $u(t,x)$ de l'\'equation de Hamilton-Jacobi ait, \`a chaque instant, un maximum \'egal \`a 0.

Nous avons alors le
\begin{theorem}
\label{t0.1}
On choisit  $R\in C^2$, et on suppose l'existence de $ I_M>0$  tel que
$\di\max_{x \in \R^d} R(x,I_M) = 0 = R(0,I_M)$. De plus, $R$ est suppos\'ee strictement concave et, pour $\vert x\vert$ grand, comprise entre deux paraboles. La donnée initiale $u_0$ est \'egalement \`a décroissance quadratique et strictement concave. 

Le probl\`eme (\ref{e0.1}) a une  unique solution $(u,I)$, o\`u $u$ est une solution classique de l'\'equation de Hamilton-Jacobi. De plus $u \in L^\infty_{\mathrm loc} \big(R^+; W^{3,\infty}_{\mathrm loc}(\R^d) \big)\cap W^{1,\infty}_{\mathrm loc} \big(R^+; L^{\infty}_{\mathrm loc}(\R^d) \big) \times W^{1,\infty}(\R)$.
\end{theorem}
L'existence pour (\ref{e0.1}) a \'et\'e d\'emontr\'ee en plusieurs endroits, voir par exemple \cite{BP} ou \cite{GB.SM.BP:09}. L'unicit\'e est donc notre r\'esultat principal, c'\'etait un probl\`eme ouvert. L'unicit\'e \'etait en effet connue seulement pour un cas tr\`es particulier (voir \cite{BP}).

Le mod\`ele (\ref{e0.1}) intervient dans la limite $\varepsilon\to 0$ des solutions de
\begin{equation} \label{e0.2}
\p_t n_\e - \e \Delta n_\e=  \frac{n_\e}{\e} R \big(x, I_\e(t) \big)\ (t>0, \; x \in \R^d),\  \       \     I_\e(t) = \int_{\R^d}  \psi(x) n_\e(t,x) dx,
\end{equation}
o\`u $n_\varepsilon(t,x)$ est la densit\'e d'une population caract\'eris\'ee par un trait biologique $x$ $d$-dimensionnel. La comp\'etition pour une ressource unique est repr\'esent\'ee par $I_\e(t)$, $\psi>0$ r\'eguli\`ere donn\'ee. Le terme $R(x,I)$ est le taux de reproduction. Les hypoth\`ses de concavit\' e sont des hypoth\`eses techniques, mais pertinentes au plan biologique. Par une transformation de Hopf-Cole $n_\e =\exp \left( {u_\e}/ {\e} \right)$ on se ram\`ene \`a l'\'equation sur $u_\e$ suivante:
\begin{equation}
\label{e0.3}
\partial_t u_\e =\e\Delta u_\e +| \nabla  u_\e |^2 +R(x,I_\e)
\end{equation}
qui, dans la limite $\varepsilon\to 0$, donne l'\'equation sur $u$. On s'attend alors \`a ce que $n_\e$ se concentre aux points o\`u $u$ est proche de 0. Et, dans cette limite, $I_\e$ appara\^{\i}t comme une sorte de multiplicateur de Lagrange.

La convergence de (\ref{e0.3}) vers (\ref{e0.1})  \`a une sous-suite pr\`es est connue depuis \cite{BP}. Le Th\'eor\`eme \ref{t0.1} donne la convergence de toute la famille, ainsi que des estimations d'erreur.  Soit $x_\e(t)$ le point o\`u $u_\e(t,.)$ atteint son maximum. On suppose l'existence de $I^0$ tel que
$
0<I^0\leq I_\e(0) := \di\int_{\R^d} \psi(x) n_\e^0(x) dx< I_M,$
et on suppose
$$
n_\e^0= e^{{u_\e^0}/{\e}}=\frac{r}{ \e^{d/2}}\, e^{{u_0}/{\e}},  \qquad \hbox{with }\;  u_0 \in C^2(\R^d) \quad \hbox{and  } \quad  \max_{x\in \R^d} u_0(x)=0.
$$
Le r\'esultat est alors le
\begin{theorem}
\label{t0.2}
Soit  $n_\e$ la solution de (\ref{e0.2}) et $u_\e$ d\'efinie par (\ref{e0.3}).  
Nous avons les d\'eveloppements asymptotiques suivants
$$
I_\e  =I+ \e I_1+ o(1),\quad x_\e =\overline x +\e \,\overline x_1+o(1), \quad u_\e=u+ \e \log (\frac{r}{ \e^{\frac d 2}})+\e\, u_1+ o(1).
$$
\end{theorem}
Les termes $I_1$, $\overline x_1$ et $u_1$ vont \^etre pr\'esent\'es dans \cite{MR2}. Ce r\'esultat implique le corollaire
\begin{corollary}
Nous avons l'approximation suivante pour $n_\e$:
$$
n_\e(t,x) = \frac{r}{ \e^{\frac d 2}} \left( \exp(u_1 +\frac {u}{\e}) +o(1) \right).
$$
En particular, lorsque $\e\to 0$, toute la suite $(n_\e)_\e$ converge:
$$
n_\e(t,x)\longrightarrow \bar \rho(t) \; \delta \big(x-\bar x (t)\big)\qquad \hbox{au sens des mesures},
$$
avec
$\rho(t)=\di \frac{I(t) }{\psi(\overline x(t))}.$
\end{corollary}

En d'autres termes, la population se concentre sur un trait dominant qui \'evolue avec le temps. On note que la convergence de $n_\e$ \`a une sous-suite pr\'es \'etait d\'ej\`a \'etablie dans \cite{LMP}.\\
 Tous ces r\'esultats seront d\'etaill\'es dans \cite{MR1} et \cite{MR2}.

\selectlanguage{english}
\section{Introduction}
The purpose of this note is to discuss uniqueness in the following problem, with unknowns $(I(t),u(t,x))$:
\begin{equation}
\label{HJ}
\left\{
\begin{array}{rll}
\partial_t u &= | \nabla  u |^2 +R(x,I)\ \  (t>0,x\in\R^d), \  \ \max_{x} u(t,x) =0,\\
I(0)&=I_0>0,\  \  u(0,x)=u_0(x),
\end{array}
\right.
\end{equation}
where $I_0>0$ and $u_0$ is a concave, quadratic function: 
$$
-\underline{L}_0 -\underline{L}_1 |x|^2 \leq u_0(x) \leq  \overline{L}_0 - \overline{L}_1 |x|^2,\       \
-2\underline{L}_1 \leq D^2u_0  \leq - 2\overline{L}_1,\       \   D^3u_0 \in L^\infty(\R^d).
$$
The constraint on the maximum of $u(t,.)$ makes the problem nonstandard. Our main result is 
\begin{theorem}
\label{thm:uniq}
Choose $R\in C^2$, and suppose  that  there is $ I_M>0$  such that
$\di\max_{x \in \R^d} R(x,I_M) = 0 = R(0,I_M)$. Also assume the following concavity and regularity properties for $R$:
$$
\begin{array}{rll}
-\underline{K}_1 |x|^2 \leq R(x,I) \leq &\overline{K}_0 -\overline{K}_1 |x|^2,  \qquad \hbox{for }\;  0\leq I \leq I_M,\\
- 2\underline{K}_1 \leq D^2 R(x,I) \leq &- 2\overline{K}_1 < 0 \   \hbox{and }D^3R(\cdot,I) \in  L^\infty(\R^d)\, \hbox{ for }0\leq I \leq I_M,\\
- \underline{K}_2\leq &\partial_I R\leq - \overline{K}_2,\\
|\partial^2_{Ix_i} R(x,I)| + | {\partial^3 R}_{Ix_{i}x_j}(x,I)| \leq &K_3,  \qquad \hbox{for  $0\leq I \leq I_M$, and $i,\,j=1,2,\cdots,d$}.
\end{array}
$$
Problem (\ref{HJ}) has a unique solution $(u,I)$, where $u$ solves the Hamilton-Jacobi equation in the classical sense. Moreover $u \in L^\infty_{\mathrm loc} \big(R^+; W^{3,\infty}_{\mathrm loc}(\R^d) \big)\cap W^{1,\infty}_{\mathrm loc} \big(R^+; L^{\infty}_{\mathrm loc}(\R^d) \big) \times W^{1,\infty}(\R)$.
\end{theorem}
Existence to (\ref{HJ}) has been proved in various contexts (see  \cite{BP,GB.SM.BP:09,LMP}). Thus, our contribution is  uniqueness, which has up to now been an open problem. The uniqueness has indeed been known only for a very particular case (see \cite{BP}).

\noindent The rest of the note is organized as follows. In Section 2, we explain the motivation and, in particular, the meaning of the various assumptions. In Section 3, we revisit existence for (\ref{HJ}), which will entail an unconventional ODE formulation for uniqueness. Section 4, which is the main part of the note, provides a fairly complete sketch of the uniqueness proof. In Section 5, we give an application.
\section{Background and motivation}
Model (\ref{HJ}) arises in the limit $\varepsilon\to 0$ of the solutions to the problem
\begin{equation} \label{para}
\p_t n_\e - \e \Delta n_\e=  \frac{n_\e}{\e} R \big(x, I_\e(t) \big)\ (t>0, \; x \in \R^d),\  \       \     I_\e(t) = \int_{\R^d}  \psi(x) n_\e(t,x) dx,
\end{equation}
where $n_\e(t,x)$ is the density of a population characterized by a $d$-dimensional biological trait $x$. The population competes for a single resource, this is represented by $I_\e(t)$,
where $\psi$ is a given positive smooth function. The term $R(x,I)$ is the reproduction rate; it is, as can be expected, very negative for large $x$ and decreases as the competition increases. Such models can be derived from individual based stochastic processes in the limit of large populations  (see \cite{NC.RF.SM:08}). The concavity assumption on $R$ is a technical one, although biologically relevant.
The Hopf-Cole transformation $n_\e =\exp \left( {u_\e}/ {\e} \right)$ yields the equation 
\begin{equation}
\label{HJE}
\partial_t u_\e =\e\Delta u_\e +| \nabla  u_\e |^2 +R(x,I_\e)
\end{equation}
 which, in the  limit $\varepsilon\to 0$, yields the equation for $u$. Now, $I_\e$ being uniformly positive and bounded in $\e$, the Hopf-Cole transformation leads to the constraint on $u$.
 Moreover, one expects that $n_\e$ concentrates at the points where $u$ is close to 0 and the function $I_\e$ appears, in the limit, as a sort of Lagrange multiplier. 
 
 This approach, based on the Hopf-Cole transformation, to study (\ref{para}) has been introduced in \cite{OD.PJ.SM.BP:05} and then developed in different contexts (see for instance \cite{BP,GB.SM.BP:09,NC.PJ:10,LMP}). Long time asymptotics of such models have also been studied in \cite{Raoulphd} and the references therein.

\section{Existence} Existence to a solution to (\ref{HJ}) is obtained by letting $\e\to 0$ in (\ref{HJE}). The main step is the
\begin{theorem}
\label{t3.1} (uniform estimates for $u_\e$, \cite{LMP})
There exists $I_m>0$ such that $0< I_m \leq I_\e (t) \leq I_M+C\e^2$. Moreover 
we have the following estimates on $u_\e$
\begin{equation}
  \label{e3.2}
  \begin{cases}
 
-\underline{L}_0 -\underline{L}_1 |x|^2  -\e 2d\underline{L}_1t \leq  u_\e(t,x) \leq \overline{L}_0 -\overline{L}_1 |x|^2 +\left( \overline{K}_0 + 2 d \e \overline{L}_1\right)t,
\\
 \underline L_1-2t\underline{K}_1 \leq D^2u_\e(t,x)  \leq -2\overline{L}_1,    \   \  \
\|D^3 u_\e(t,\cdot) \|_{L^\infty} \leq C(T), \qquad \hbox{for $t\in [0,T]$}. 
\end{cases}
\end{equation}
\end{theorem}
The bounds for $u_\e$ can be obtained for any uniformly bounded function $I_\e$,
not only for that of (\ref{para}). This remark will be an important ingredient of the uniqueness proof.  

\section{Uniqueness} For a given continuous function $I(t)$ such that $0<I(t)<I_M$, one may construct a solution of $\partial_t u =| \nabla  u |^2 +R(x,I)$ with initial datum $u_0$.  Just as in Theorem \ref{t3.1}, this solution satisfies estimates (\ref{e3.2}). And so, $u(t,.)$ being strictly concave and quadratically decreasing, there exists a unique function $\bar x(t)$ such that $u(t,\bar x(t))=\di\max_{x\in\RR^d}u(t,x)$.  Assume that $I(t)$ is chosen such that $u(t,\bar x(t))=0$. Then, from the equation on $u$ we deduce that $R(\bar x(t),I(t))~=0$. Notice also that, because $\partial_I R<0$, we have $R(\bar x(t),0)>0$. 
 Finally, differentiating $\nabla u(t,\bar x(t))=0$ and plugging in the equation for $u$ we obtain an ODE for $\bar x$: $\dot{\bar x}(t) = \left( -D^2u \big(t, \bar x(t)\big) \right)^{-1} \nabla_x R\big(\bar x(t),\bar I(t) \big)$. 

The idea is thus to change the constrained problem (\ref{HJ}) by the following slightly nonstandard differential system:
\begin{equation}
\label{system}
\begin{cases}
R \left(\overline x(t),I(t) \right)=0,& \hbox{for $t\in [0,T]$},\\
\dot{\overline x}(t) = \left( -D^2u \big(t, \bar x(t)\big) \right)^{-1} \nabla_x R\big(\bar x(t),\bar I(t) \big),& \hbox{for $t\in [0,T]$},\\
\partial_t u = | \nabla  u |^2 +R(x,I),& \hbox{in $ [0,T] \times \R^d$},
\end{cases}
\end{equation}
with initial conditions
\begin{equation}
\label{first}
I(0)=I_0,\qquad u(0,\cdot)=u_0(\cdot),\qquad \overline x(0)=\overline x_0,\qquad \hbox{such that $R(x_0,I_0)=0$.} 
\end{equation}
Note that (\ref{system}) is really a differential system because the assumptions on $R$ imply that $I(t)$ can implicitely be expressed in terms of $\bar x(t)$. And it is slightly nonstandard because $\bar x$ solves an ODE whose nonlinearity depends on $u$. Finally, note that, as soon as $u$ satisfies the concavity and regularity estimates (\ref{e3.2}), system (\ref{system}) is equivalent to the constrained problem (\ref{HJ}).

This suggests to use a simple fixed point argument to prove uniqueness to  (\ref{system})  (and so, to (\ref{HJ})). Which in turn suggests to set up the following scheme: starting from $ x(t)\in C([0,T],\RR^d)$, such that $x(0)=x_0$, where $R(x_0,0)>0$. Let $I(t)$ solve $R(x(t),I(t))=0$ on $[0,T]$ with $R(x_0,I_0)=0$. Let $v(t,.)$ be the unique solution to 
\begin{equation}
\label{e4.3}
\partial_t v = | \nabla  v |^2 +R(x,I),\     \   v(0,x)=u_0(x).
\end{equation}
Let $y(t)$ solve  $\dot{y}(t) = \left( -D^2v \big(t, x(t)\big) \right)^{-1} \nabla_x R\big(x(t),I(t) \big)$ on $[0,T]$ with initial datum $y(0)=x_0$.  Setting $y:=\Phi(x)$, we notice that uniqueness is proved as soon as we have proved that $\Phi$ has a unique fixed point. One additional feature about a solution $(\bar I,u,\bar x)$ of (\ref{system}):
\begin{lemma}
\label{l4.1}
The function $\bar I(t)$ is increasing.
\end{lemma}
We claim that our problem reduces to proving the 
\begin{theorem}
\label{t4.1}
There exists  $C>0$ universal and $\delta>0$, which is small as $R(x_0,0)$ tends to 0, such that $\Phi$ is a contraction from $C([0,\delta],\overline B_{C}(x_0))$ to itself;
here $B_r(a)$ denotes the ball of centre $a\in\RR^d$ and radius $r>0$. 
\end{theorem}
Note indeed that,   by Lemma \ref{l4.1},  we have, because $\partial_IR<0$:
$$
R(\bar x(\delta),0)=R(\bar x(\delta),0)-R(\bar x(\delta),\bar I(\delta))\geq c\bar I(\delta)\geq cI_0,
$$
for some universal $c>0$. Hence Theorem \ref{t4.1} can be iterated to yield global existence and uniqueness.

Let us give an overview of the proof of Theorem \ref{t4.1}. For $I\in \mathrm{C} \left( [0,\delta] ; [0,I_M] \right)$, let $V(I)$ be the (unique) solution of (\ref{e4.3}). The main step is the 
following
\begin{lemma}
\label{l4.2}
Let $I_1,\, I_2 \in \mathrm{C} \left( [0,\delta] ; [0,I_M] \right)$. Then
$$
\|V(I_1)-V(I_2)\|_{W^{2,\infty}([0,\delta]\times \R^d)} \leq C \| I_1 -I_2 \|_{L^\infty([0,\delta])} \delta.
$$
\end{lemma}
This lemma, once proved, opens the way to Theorem \ref{t4.1}. Indeed the equation $R(x,I)=0$ yields a smooth mapping $x\mapsto I$, and $I \mapsto V$ is a Lipschitz mapping thanks to  Lemma \ref{l4.2}. Moreover, the equation $\dot{y}(t) = \left( -D^2v \big(t, x(t)\big) \right)^{-1} \nabla_x R\big(x(t),I(t) \big)$ yields a Lipschitz mapping $v\mapsto y$ by the estimates for $v$ given by Theorem \ref{t3.1}.

Lemma \ref{l4.2} is more involved. If $I_1$ and $I_2$ are as in the assumptions of the lemma, the function $r=V(I_1)-V(I_2)$ solves 
\begin{equation}
\label{e4.4}
\begin{cases}
\partial_t r= \left( \nabla v_1 + \nabla v_2 \right) \cdot \nabla r +R(x,I_2) -R(x,I_1),& \hbox{in $[0,\delta]\times \R^d$}\\
r(0,x)=0,& \hbox{for all $x\in \R^d$.}
\end{cases}
\end{equation}
with $v_i=V(I_i)$. Note that the above equation has a unique classical solution which can be computed by the method of characteristics. The characteristics solve
\begin{equation}
\label{e4.5}
\dot{\gamma }(t) = - \nabla v_1 (t,\gamma)-\nabla v_2   (t,\gamma),
\end{equation}
and, due to the estimates of Theorem \ref{t3.1}, they exist globally. So, one may successively express $r$ given by integration along characteristics, and estimate its derivatives recursively.
\section{Application}
Convergence of  (\ref{HJE}) to (\ref{HJ}) had already been proved in \cite{BP,GB.SM.BP:09}, along subsequences. The uniqueness part of Theorem \ref{thm:uniq} yields the convergence of the full family of  solutions $u_\e$ of (\ref{HJE}), instead of convergence along a subsequence. Moreover it allows an expansion of $I_\e$, $u_\e$ and $x_\e$ (the maximum point of $u_\e$ at each time) in terms of $\e$.  Here are the results.

Assume that there is  $I^0$ such that
$
0<I^0\leq I_\e(0) := \di\int_{\R^d} \psi(x) n_\e^0(x) dx< I_M,$
and that
$$
n_\e^0= e^{{u_\e^0}/{\e}}=\frac{r}{ \e^{d/2}}\, e^{{u^0}/{\e}},  \qquad \hbox{with }\;  u^0 \in C^2(\R^d) \quad \hbox{and  } \quad  \max_{x\in \R^d} u^0(x)=0.
$$
The result is the
\begin{theorem}
\label{thm:approx}
Let $n_\e$ be the solution of (\ref{para}) and $u_\e$ be defined by (\ref{HJE}).  
We have the following asymptotic expansions
$$
I_\e  =I+ \e I_1+ o(1),\quad x_\e =\overline x +\e \,\overline x_1+o(1), \quad u_\e=u+ \e \log (\frac{r}{ \e^{\frac d 2}})+\e\, u_1+ o(1).
$$
\end{theorem}
The terms $I_1$, $\overline x_1$ and $u_1$ will be provided in \cite{MR2}. This yields the corollary
\begin{corollary}
We have the following approximation for $n_\e$:
$$
n_\e(t,x) = \frac{r}{ \e^{\frac d 2}} \Big( \exp\big(u_1(t,x) +\frac {u(t,x)}{\e}\big) +o(1) \Big).
$$
In particular, as $\e\to 0$, the whole sequence $(n_\e)_\e$ converges:
$$
n_\e(t,x)\longrightarrow \bar \rho(t) \; \delta \big(x-\bar x (t)\big),\qquad \hbox{weakly in the sense of measures},
$$
with 
$\rho(t)=\di \frac{I(t) }{\psi(\overline x(t))}.$
\end{corollary}

In other words, the population density concentrates on a dominant trait which evolves in time. We note that the convergence of $n_\e$ along subsequences was already established in \cite{LMP}.\\
 The above results will be detailed in \cite{MR1} and \cite{MR2}.
 
\section*{Acknowledgements}
S. Mirrahimi was partially funded by the french ANR projects KIBORD ANR-13-BS01-0004 and MODEVOL ANR-13-JS01-0009.
J.-M. Roquejoffre  has received funding from the European Research Council
under the European Union's Seventh Framework Programme (FP/2007-2013) / ERC Grant
Agreement n. 321186 - ReaDi - ``Reaction-Diffusion Equations, Propagation and Modelling'' held by Henri Berestycki.
He was also partially supported by the French National Research Agency
(ANR), within the project NONLOCAL ANR-14-CE25-0013.  Both authors are grateful to the Labex CIMI for organizing a
PDE-probability quarter in Toulouse, which provided quite a stimulating environment to this collaboration.


\begin{thebibliography}{00}



\bibitem{GB.SM.BP:09}
{\sc G.~Barles, S.~Mirrahimi, and B.~Perthame},
{\it Concentration in {L}otka-{V}olterra parabolic or integral equations:
  a general convergence result.}
 Methods Appl. Anal., {\bf 16(3)} (2009) pp.321--340.
 


\bibitem{NC.RF.SM:08}
{\sc N.~Champagnat, R.~Ferri\`ere, and S.~M\'el\'eard},
 {\it Individual-based probabilistic models of adaptive evolution and
  various scaling approximations}, {\bf 59} (2008) Progress in Probability. Birkha\"{u}ser.

\bibitem{NC.PJ:10}
{\sc N.~Champagnat and P.-E. Jabin},
{\it The evolutionary limit for models of populations interacting competitively via several resources,}
Journal of Differential Equations, {\bf 261} (2011), pp.179--195.
%


\bibitem{OD.PJ.SM.BP:05}
{\sc O.~Diekmann, P.-E. Jabin, S.~Mischler, and B.~Perthame},
{\it The dynamics of adaptation: an illuminating example and a
  {H}amilton-{J}acobi approach},
 Th. Pop. Biol., {\bf 67(4)} (2005) pp.257--271.
 
 
\bibitem{LMP} {\sc A. Lorz, S. Mirrahimi, B Perthame}, {\it Dirac mass dynamics in a multidimensional nonlocal parabolic equation.} Communications in Partial Differential Equations, {\bf 36} (2011) pp.1071--1098.
\bibitem{MR1} {\sc S. Mirrahimi, J.-M. Roquejoffre}, {\it Uniqueness in a Hamilton-Jacobi equation with constraints,} in preparation.
\bibitem{MR2}  {\sc S. Mirrahimi, J.-M. Roquejoffre}, {\it Approximation of solutions of selection-mutation models and error estimates,} in preparation.
\bibitem{BP} {\sc  B. Perthame, G. Barles}, {\it Dirac concentrations in Lotka-Volterra parabolic PDEs,} Indiana Univ. Math. J. {\bf 57} (2008), pp. 3275--3301.


\bibitem{Raoulphd}
{\sc G.~Raoul},
{\it Etude qualitative et num\'{e}rique d'\'{e}quations aux
  d\'{e}riv\'{e}es partielles issues des sciences de la nature},
PhD thesis, ENS Cachan, 2009.



\end{thebibliography}
\end{document}